\documentclass[12pt]{amsart}
\usepackage{amssymb,amsmath,color}

\theoremstyle{plain}
\newtheorem{thm}{Theorem}[section]

\theoremstyle{definition}
\newtheorem{ex}{Example}[section]

\renewcommand{\leq}{\leqslant}
\renewcommand{\geq}{\geqslant}

\numberwithin{equation}{section}

\begin{document}
\title[The spectral density function of the Bochner Laplacian]{The spectral density function \\ of the renormalized Bochner Laplacian\\ on a symplectic manifold}

\author[Y. A. Kordyukov]{Yuri A. Kordyukov}
\address{Institute of Mathematics, Ufa Federal Research Centre, Russian Academy of Sciences, 112~Chernyshevsky St., 450008 Ufa, Russia and Novosibirsk State University, 1, Pirogova St, 630090, Novosibirsk, Russia} \email{yurikor@matem.anrb.ru}

\thanks{}



\begin{abstract}
We consider the renormalized Bochner Laplacian acting on tensor powers of a positive line bundle on a compact symplectic manifold. We derive an explicit local formula for the spectral density function in terms of coefficients of the Riemannian metric and symplectic form.
\end{abstract}

\allowdisplaybreaks

\date{}

 \maketitle

\section{Introduction}
\subsection{Preliminaries}
 Let $(X,\mathbf B)$ be a compact symplectic manifold of dimension $2n$. We assume that there exists a Hermitian line bundle $(L,h^L)$ on $X$ with a Hermitian connection $\nabla^L : C^\infty(X,L)\to C^\infty(X,T^*X\otimes L)$, which satisfies the pre-quantization condition:
\begin{equation}\label{e:def-omega}
iR^L=\mathbf B,
\end{equation}
where $R^L=(\nabla^L)^2$ is the curvature of the connection. Thus, $[\mathbf B]\in H^2(X,2\pi \mathbb Z)$.

Let $g$ be a Riemannian metric on $X$. For any $p\in \mathbb N$, denote by $L^p:=L^{\otimes p}$ be the $p$-th tensor power of $L$. Consider the induced Bochner Laplacian $\Delta^{L^p}$ acting on $C^\infty(X,L^p)$ by
\begin{equation}\label{e:DeltaLp}
\Delta^{L^p}=\big(\nabla^{L^p}\big)^{\!*}\,
\nabla^{L^p},
\end{equation}
where $\nabla^{L^p}: {C}^\infty(X,L^p)\to {C}^\infty(X, T^*X \otimes L^p)$ is the connection on $L^p$ induced by $\nabla^{L}$, and $(\nabla^{L^p})^{*}:
{C}^\infty(X,T^*X\otimes L^p)\to {C}^\infty(X,L^p)$ is the formal adjoint of $\nabla^{L^p}$.  

In the case where $(L,h^L)$ is the trivial Hermitian line bundle, the Hermitian connection $\nabla^L$ can be written as $\nabla^L=d-i \mathbf A$ with some real-valued 1-form $\mathbf A$, and we have
\[
R^L=-id\mathbf A,\quad \mathbf B=d\mathbf A. 
\]
Thus, $\mathbf B$ can be considered as a magnetic 2-form and $\mathbf A$ as the associated magnetic potential. The Bochner Laplacian $\Delta^{L^p}$ is related with the semiclassical magnetic Schr\"odinger operator
\[
H^{\hbar}_{\mathbf A}=(i\hbar d+\mathbf A)^*(i\hbar d+\mathbf A), \quad \hbar>0
\]
by the formula
\[
\Delta^{L^p}=\hbar^{-2}H^{\hbar}_{\mathbf A}, \quad \hbar=\frac{1}{p},\quad p\in \mathbb N. 
\]
 
Let $B\in \operatorname{End}(TX)$ be a skew-adjoint endomorphism such that 
\begin{equation}\label{e:defJ0}
\mathbf B(u,v)=g(Bu,v), \quad u,v\in TX. 
\end{equation}

The renormalized Bochner Laplacian $\Delta_p$ is a second order differential operator acting on $C^\infty(X,L^p)$ by
 \[
\Delta_p=\Delta^{L^p}-p\tau,
 \]
where $\tau$ is a smooth function on $X$ given by 
 \begin{equation}\label{e:def-tau}
 \tau(x)=\frac 12 \operatorname{Tr}[(B(x)^*B(x))^{1/2}],\quad x\in X.
 \end{equation}
This operator was introduced in \cite{Gu-Uribe}. 
 
An almost complex structure $J\in \operatorname{End}(TX)$ compatible with $\mathbf B$ and $g$ is defined by 
\begin{equation}\label{e:defJ}
J=B(B^*B)^{-1/2}.
\end{equation} 
We put
  \begin{equation}\label{e:def-mu0}
 \mu_0=\inf_{u\in T_xX, x\in X}\frac{\mathbf B_x(u,J(x)u)}{|u|_g^2}.
 \end{equation}
We denote by $\sigma(\Delta_p)$ the spectrum of $\Delta_p$ in $L^2(X,L^p)$. By \cite[Corollary 1.2]{ma-ma02}, there exists a constant $C_L>0$ such that for any $p$
 \[
 \sigma(\Delta_p)\subset [-C_L,C_L]\cup [2p\mu_0-C_L,+\infty).
 \]
Consider the finite-dimensional vector subspace $\mathcal H_p\subset L^2(X,L^p)$ spanned by the eigensections of $\Delta_p$ corresponding to eigenvalues in $[-C_L,C_L]$. Its dimension $d_p$ grows polynomially as $p\to \infty$ (\cite[Corollary 1.2]{ma-ma02}).
Denote by $\lambda_j^{(p)}, j=1,2,\ldots,d_p$, the eigenvalues of $\Delta_p$ in $[-C_L,C_L]$ taken with multiplicities. The spectral density function is a function $\rho\in C^\infty(X)$ such that for any $f\in C(\mathbb R)$
\[
\lim_{p\to \infty}\frac{1}{d_p}\sum_{j=1}^{d_p}f(\lambda_j^{(p)})=\frac{1}{{\rm vol}(X)}\int_X f(\rho(x))d\mu_X(x),
\]
where $d\mu_X=\frac{1}{n!}\mathbf B^n$ is the Liouville volume form associated with the symplectic form $\mathbf B$ and ${\rm vol}(X)=\int_Xd\mu_X$. The existence and uniqueness of such a function $\rho$ was shown in \cite{Gu-Uribe}. 

The main result of the paper is an explicit local formula for the function $\rho$ in terms of the coefficients of the tensors $g$ and $\mathbf B$. 

\subsection{The main result}
We consider the complexified tangent bundle $TX\otimes\mathbb C$ equipped with the $\mathbb C$-bilinear form $\langle\cdot,\cdot\rangle$ induced by $g$.  Define  $\mathcal J\in \operatorname{End}(TX\otimes\mathbb C)$ by 
\begin{equation}\label{e:defcJ0}
\mathcal J=-iB.
\end{equation}
It is skew-adjoint with respect to $\langle\cdot,\cdot\rangle$.

Fix an arbitrary point $x_0\in X$. The almost complex structure $J_{x_0}$ defined in \eqref{e:defJ} induces a splitting $T_{x_0}X\otimes \mathbb C=T^{(1,0)}_{x_0}X\oplus T^{(0,1)}_{x_0}X$, where $T^{(1,0)}_{x_0}X$ and $T^{(0,1)}_{x_0}X$ are the eigenspaces of $J_{x_0}$ corresponding to its eigenvalues $i$ and $-i$ respectively. The operator $\mathcal J_{x_0}$ gives rise to a positive operator $\mathcal J_{x_0} : T^{(1,0)}_{x_0}X\to T^{(1,0)}_{x_0}X$. We choose an orthonormal base $\{w_j : j=1,\ldots,n\}$ of $T^{(1,0)}_{x_0}X$, consisting of its eigenvectors:
\[
\mathcal J_{x_0}w_j=a_jw_j,\quad a_j>0, \quad j=1,\ldots,n.
\]
Then the vectors  $\{e_k : k=1,\ldots,2n\}$ given by
\begin{equation}\label{e:ej}
e_{2j-1}=\frac{1}{\sqrt{2}}(w_j+\bar w_j), \quad e_{2j}=\frac{i}{\sqrt{2}}(w_j-\bar w_j), \quad j=1,\ldots,n,
\end{equation}
form an orthonormal base of $T_{x_0}X$. By means of this base, we identify $T_{x_0}X$ with $\mathbb R^{2n}$. We use the corresponding coordinates on $T_{x_0}X$: $Z=(Z_1,\ldots,Z_{2n})\in \mathbb R^{2n}\mapsto \sum_{j=1}^{2n}Z_je_j\in T_{x_0}X$, as well as the complex coordinates $z\in\mathbb C^{n}$, $z_j=Z_{2j-1}+iZ_{2j}, j=1,\ldots,n$. We put
\[
\frac{\partial}{\partial z_j}=\frac{1}{2}\left(\frac{\partial}{\partial Z_{2j-1}}-i\frac{\partial}{\partial Z_{2j}}\right), \quad \frac{\partial}{\partial \bar{z}_j}=\frac{1}{2}\left(\frac{\partial}{\partial Z_{2j-1}}+i\frac{\partial}{\partial Z_{2j}}\right).
\]
Thus, $w_j=\sqrt{2}\tfrac{\partial}{\partial z_j}$, $\bar w_j=\sqrt{2}\tfrac{\partial}{\partial\bar z_j}$. 

Consider the endomorphism $|\mathcal J|:TX\otimes\mathbb C \to TX\otimes\mathbb C$ given by
\[
|\mathcal J|:=(\mathcal J^*\mathcal J)^{1/2}=-iJ\mathcal J.
\]
We denote by $\nabla^{TX}$ the Levi-Civita connection on $(X,g)$. For each tensor $\psi$ on $X$, we denote by $\nabla^X\psi$ the covariant derivative of $\psi$ induced by $\nabla^{TX}$. Thus, $\nabla^X\mathcal J, \nabla^XJ\in T^*X\otimes \operatorname{End}(TX)$, $\nabla^X\nabla^X\mathcal J, \nabla^XJ\in T^*X\otimes T^*X\otimes \operatorname{End}(TX)$. 

\begin{thm}\label{t:main}
The following formula holds:
\begin{equation}\label{e:rho}
\begin{aligned}
\rho(x_0)=  & -\sum_{j,k=1}^n\frac{8}{a_j+a_k} \left \langle (\nabla ^{X} \nabla ^{X}\mathcal J)_{(\tfrac{\partial}{\partial z_j},\tfrac{\partial}{\partial z_k})} \tfrac{\partial}{\partial \bar{z}_j},  \tfrac{\partial}{\partial \bar{z}_k} \right \rangle \\ 
 & + \sum_{j=1}^n \frac{1}{a_j} \left( {\mbox{\rm tr}}_{|TX} (\nabla ^{X} \nabla ^{X}\mathcal J)_{(\tfrac{\partial}{\partial z_j},\tfrac{\partial}{\partial \bar z_j})}
 - {\mbox{\rm tr}}_{|TX} \Big(\nabla ^{X} \nabla ^{X}|\mathcal J|\Big)_{(\tfrac{\partial}{\partial z_j},\tfrac{\partial}{\partial \bar z_j})}\right). \\
&  +\sum_{j,k,\ell=1}^n \frac{8}{9a_k(a_j+a_\ell)}\left|\left\langle  (\nabla ^{X}_{\tfrac{\partial}{\partial z_k}}\mathcal J )\tfrac{\partial}{\partial \bar z_j},\tfrac{\partial} {\partial \bar z_\ell}\right\rangle\right|^2 \\ & + \sum_{j, k,\ell=1}^n \frac{8}{a_k(a_j+a_k+a_\ell)}\left|\left\langle(\nabla ^{X}_{\tfrac{\partial}{\partial \bar z_k}}\mathcal J)\tfrac{\partial}{\partial \bar z_j},\tfrac{\partial} {\partial\bar{z}_\ell}\right\rangle \right|^2,
\end{aligned}
\end{equation}
where all the tensors are evaluated at $x_0$ and the covariant derivatives are taken with respect to normal coordinates at $x_0$ defined by the base \eqref{e:ej}.
\end{thm}

Let us mention two particular cases of the main result. 

\begin{ex}\label{ex:1} Assume that $B=2\pi J$ (the almost-K\"ahler case). Then $|\mathcal J|=2\pi I$, $a_j=2\pi$ and we get the formula proved in \cite[formula (2.30)]{ma-ma08a}:
\[
\rho(x_0) =\frac{1}{24}|\nabla^XJ|^2.
\]
\end{ex}

\begin{ex}\label{ex:2}
Assume that $\nabla^X J=0$. For instance, if, in addition, $J$ is a complex structure, then $(X,J,g)$ is K\"ahler. Then formula \eqref{e:rho} takes the form:
\begin{multline*}
\rho(x_0) = 
 \sum_{j,k=1}^n\frac{8}{a_j+a_k}\left \langle (\nabla ^{X} \nabla ^{X} |\mathcal J| )_{(\tfrac{\partial}{\partial z_j},\tfrac{\partial}{\partial \bar z_j})}\tfrac{\partial}{\partial \bar z_k} -(\nabla ^{X} \nabla ^{X} |\mathcal J|)_{(\tfrac{\partial}{\partial z_j},\tfrac{\partial}{\partial \bar z_k})}\tfrac{\partial}{\partial \bar z_j}, \tfrac{\partial}{\partial z_k} \right \rangle\\  +\sum_{j,k,\ell=1}^n \frac{8}{9a_k(a_j+a_\ell)}\left|\left\langle  \nabla ^{X}_{\tfrac{\partial}{\partial \bar z_\ell}}|\mathcal J|\tfrac{\partial} {\partial\bar{z}_j}- \nabla ^{X}_{\tfrac{\partial}{\partial \bar z_j}}|\mathcal J|\tfrac{\partial} {\partial\bar{z}_\ell}, \tfrac{\partial}{\partial z_k}\right\rangle\right|^2.
\end{multline*}
\end{ex}

The spectral density function $\rho$ is a spectral invariant of the renormalized Bochner Laplacian and, therefore, a differential geometric invariant of the pair $(g,\mathbf B)$ on the compact manifold $X$. There is a well-elaborated theory of differential geometric invariants of Riemannian manifolds based on the notion of Riemann curvature tensor, the theory of characterictic classes and Weyl invariance theory. The theory of spectral invariants of geometric differential operators (like the Laplace-Beltrami operator or the Hodge Laplacian on differential forms) on Riemannian manifolds is also quite well studied (see, for instance, \cite{Gilkey95} and references therein). On the other hand, differential geometric invariants of pairs $(g,\mathbf B)$, consisting of a Riemannian metric $g$ and a symplectic form $\mathbf B$ on a smooth manifold, as well as spectral invariants of geometric differential operators associated with such pairs (like the Bochner Laplacian) are much less studied (for some information, see \cite{ColinDeV12,karasev07,karasev-osborn05,KT-19} and references therein). Theorem~\ref{t:main} gives an explicit expression for a particular example of a spectral invariant associated with a pair $(g,\mathbf B)$ in terms of geometric data. It would be highly interesting to give an interpretation of this invariant in terms of differential geometric invariants of the pair $(g,\mathbf B)$.

The proof of Theorem~\ref{t:main} is based on the methods and results developed in \cite{ma-ma:book,ma-ma08a}. First, we use the fact that the spectral density function coincides with the leading coefficient in the asymptotic expansion of the generalized Bergman kernel $P_{1,p}$ of the renormalized Bochner Laplacian $\Delta_p$ as $p\to \infty$, that was first observed in \cite{Gu-Uribe} and extended and refined in \cite{ma-ma08a}.  
A method to compute the coefficients in the asymptotic expansion of the generalized Bergman kernel by recurrence was developed by Ma and Marinescu in \cite{ma-ma08a}, where an integral formula for the first two
coefficients was obtained and some coefficients were explicitly computed under the assumption $B=2\pi J$. It should be noted that the paper \cite{B-Uribe02} is also dedicated to the computation of the spectral density function under assumption $B=2\pi J$ (concerning this computation, see \cite[Remark 3.2]{ma-ma08a}). 

In this paper, we start with the formula due to Ma-Marinescu and complete the computation of the leading coefficient in the general case, which results in Theorem~\ref{t:main}.

\section{Preliminaries}
\subsection{The model operator}
Let $x_0\in X$ be an arbitrary point, which will be fixed from now on. We consider the second order differential operator $\mathcal L_{x_0}$ in $C^\infty(T_{x_0}X)$ given by
\begin{equation}\label{e:defL0}
\mathcal L_{x_0}=-\sum_{j=1}^{2n} \left(\nabla_{e_j}-\tfrac 12i \mathbf B_{x_0}(Z,e_j)\right)^2-\tau (x_0), 
\end{equation}
where $\{e_j\}_{j=1,\ldots,2n}$ is an orthonormal base in $T_{x_0}X$. For $U\in T_{x_0}X$ we denote by $\nabla_U$ the ordinary operator of differentiation in the direction $U$ on $C^\infty(T_{x_0}X)$. The operator $\mathcal L_{x_0}$ is well-defined, i.e., it is independent of the choice of the base. 

From now on, we will assume that an orthonormal base $\{e_j\}_{j=1,\ldots,2n}$ in $T_{x_0}X$ is given by \eqref{e:ej} and use the associated coordinates $Z\in \mathbb R^{2n}$ on $T_{x_0}X$ as well as the complex coordinates $z\in\mathbb C^{n}$. Thus, we have
\begin{align}
J_{x_0}\frac{\partial}{\partial z_j} =& i\frac{\partial}{\partial z_j},\quad J_{x_0}\frac{\partial}{\partial \bar z_j}=-i\frac{\partial}{\partial \bar z_j},  \quad j=1,\ldots,n,\label{e:eigenJ}\\
\mathcal J_{x_0}\frac{\partial}{\partial z_j}=& a_j\frac{\partial}{\partial z_j},\quad \mathcal J_{x_0}\frac{\partial}{\partial \bar z_j}=-a_j\frac{\partial}{\partial \bar z_j},  \quad j=1,\ldots,n,\label{e:eigencJ}\\
|\mathcal J|\frac{\partial}{\partial z_j}=& a_j\frac{\partial}{\partial z_j}, \quad |\mathcal J|\frac{\partial}{\partial \bar z_j} =a_j\frac{\partial}{\partial \bar z_j}, \quad j=1,\ldots,n,\label{e:eigen|J|}
\end{align}
and 
\[
 \tau(x_0)=\frac 12 \operatorname{Tr}(|\mathcal J_{x_0}|)=\sum_{j=1}^na_j.
\]
As in \cite[Section 1.4]{ma-ma08a}, we write $\mathcal L=\mathcal L_{x_0}$ in terms of the creation and annihilation operators. We define first order differential operators $b_j,b^{+}_j, j=1,\ldots,n,$ on $\mathbb R^{2n}\cong T_{x_0}X$ by
\[
b_j= -2{\frac{\partial}{\partial z_j}}+\frac{1}{2}a_j\bar{z}_j,\quad
b^{+}_j=2{\frac{\partial}{\partial\bar{z}_j}}+\frac{1}{2}a_j z_j, \quad j=1,\ldots,n.
\]
Then $b^{+}_j$ is the formal adjoint of $b_j$ on $L^2({\mathbb R}^{2n})$, and ${\mathcal L}=\sum_{j=1}^n b_j b^{+}_j$. 

We define a function $\mathcal P_{x_0}\in C^\infty({\mathbb R}^{2n}\times {\mathbb R}^{2n})$ by
\begin{equation}\label{e:defP} 
\mathcal P_{x_0}(Z, Z^\prime)=\frac{1}{(2\pi)^n}\prod_{j=1}^na_j \exp\left(-\frac 14\sum_{k=1}^na_k(|z_k|^2+|z_k^\prime|^2- 2z_k\bar z_k^\prime) \right).
\end{equation}
It is the Bergman kernel of the operator $\mathcal L$, i.e., the smooth kernel with respect to $dZ$ of the orthogonal projection $\mathcal P=\mathcal P_{x_0}$ in $L^2({\mathbb R}^{2n})$ to the kernel of $\mathcal L$.

\subsection{Relation with generalized Bergman kernels}
Here, we recall a description of the spectral density function in terms of the asymptotic expansion of the generalized Bergman kernel $P_{1,p}$ of the renormalized Bochner Laplacian $\Delta_p$ as $p\to \infty$, first observed in \cite{Gu-Uribe} and extended and refined in \cite{ma-ma08a}.

First, we introduce normal coordinates near $x_0$. Let $a^X$ be the injectivity radius of $(X,g)$. We will identify $B^{T_{x_0}X}(0,a^X)$ with $B^{X}(x_0,a^X)$ by the exponential map $\operatorname{exp}^X : T_{x_0}X\to X$. For $Z\in B^{T_{x_0}X}(0,a^X)$ we identify $L_Z$ to $L_{x_0}$ by parallel transport with respect to the connection $\nabla^L$ along the curve $\gamma_Z : [0,1]\ni u \to \exp^X_{x_0}(uZ)$. We consider the line bundle $L_0$ with fibers $L_{x_0}$ on $T_{x_0}X$ and denote by $\nabla^L$, $h^L$ the connection and the metric on the restriction of $L_0$ to $B^{T_{x_0}X}(0,a^X)$ induced by the identification $B^{T_{x_0}X}(0,a^X)\cong B^{X}(x_0,a^X)$ and the trivialization of $L$ over $B^{T_{x_0}X}(0,a^X)$. 

Let $dv_{TX}$ be the Riemannian volume form of $(T_{x_0}X, g^{T_{x_0}X})$ and let $dv_X$ be the volume form on $B^{T_{x_0}X}(0,a^X)$, corresponding to the Riemannian volume form $dv_X$ on $B^{X}(x_0,a^X)$ under identification $B^{T_{x_0}X}(0,a^X)\cong B^{X}(x_0,a^X)$. Let $\kappa_{x_0}$ be the smooth positive function on $B^{T_{x_0}X}(0,a^X)$ defined by the equation
\[
dv_{X}(Z)=\kappa_{x_0}(Z)dv_{TX}(Z), \quad Z\in B^{T_{x_0}X}(0,a^X). 
\] 

Let $P_{\mathcal H_p}$ be the orthogonal projection from $L^2(X,L^p)$ onto $\mathcal H_p$.  We denote by $P_{q,p}(x,x^\prime)$, $x,x^\prime\in X$ the smooth kernel of the operator $\Delta_p^qP_{\mathcal H_p}$ with respect to the Riemannian volume form $dv_X$, which is called \emph{a generalized Bergman kernel} of $\Delta_p$. Under our trivialization, it induces a smooth function $P_{q,p,x_0}(Z,Z^\prime)$ on the set of all $Z,Z^\prime\in T_{x_0}X$ with $x_0\in X$ and $|Z|, |Z^\prime|<a_X$. For any $k\in \mathbb N$ and $x_0\in X$,   we have (cf. \cite[Theorem 1.1]{bergman})
\begin{equation}\label{e:Pqr}
p^{-n}P_{q,p,x_0}(Z,Z^\prime)\cong 
\sum_{r=2q}^kF_{q,r,x_0}(\sqrt{p} Z, \sqrt{p}Z^\prime)\kappa_{x_0}^{-\tfrac 12}(Z)\kappa_{x_0}^{-\tfrac 12}(Z^\prime)p^{-\tfrac{r}{2}+q},
\end{equation}
where
\begin{equation}\label{e:Fqr}
F_{q,r,x_0}(Z,Z^\prime)=J_{q,r,x_0}(Z,Z^\prime)\mathcal P_{x_0}(Z,Z^\prime),
\end{equation}
$J_{q,r,x_0}(Z,Z^\prime)$ are polynomials in $Z, Z^\prime$, depending smoothly on $x_0$, with the same parity as $r$ and $\operatorname{deg} J_{q,r,x_0}\leq 3r$.

This means that there exist $\varepsilon\in (0,a_X]$ and $C_0>0$ with the following property:
for any $l\in \mathbb N$, there exist $C>0$ and $M>0$ such that for any $x_0\in X$, $p\geq 1$ and $Z,Z^\prime\in T_{x_0}X$, $|Z|, |Z^\prime|<\varepsilon$, we have 
\begin{multline*}
\Bigg|p^{-n}P_{q,p,x_0}(Z,Z^\prime)\kappa_{x_0}^{\tfrac 12}(Z)\kappa_{x_0}^{\tfrac 12}(Z^\prime) -\sum_{r=2q}^kF_{q,r,x_0}(\sqrt{p} Z, \sqrt{p}Z^\prime)p^{-\tfrac{r}{2}+q}\Bigg|_{\mathcal C^{l}(X)}\\ 
\leq Cp^{-\tfrac{k+1}{2}+q}(1+\sqrt{p}|Z|+\sqrt{p}|Z^\prime|)^M\exp(-\sqrt{C_0p}|Z-Z^\prime|)+\mathcal O(p^{-\infty}).
\end{multline*}
This estimate was introduced in \cite{dai-liu-ma} for the spin$^c$ Dirac operator and $q=0$ and in \cite{ma-ma:book,ma-ma08a} for the K\"ahler case and arbitrary $q$ (cf. also \cite{LMM} for the renormalized Bochner Laplacian in the case $q=0$). 

The leading coefficients of these expansions satisfy 
\begin{equation}\label{e:Jqq}
J_{q,2q,x_0}(Z,Z^\prime)=J_{q,2q,x_0}(0,0)=(J_{1,2,x_0}(0,0))^q, \quad Z,Z^\prime\in T_{x_0}X.
\end{equation}
To prove these identities, one can use Toeplitz operator calculus introduced in \cite{ma-ma:book,ma-ma08} for spin$^c$ Dirac operator and K\"ahler case (also with an auxiliary bundle) and extended to the case under consideration in \cite{ioos-lu-ma-ma,bergman}. Using the asymptotic expansions of the generalized Bergman kernels mentioned above and the characterization of Toeplitz operators in terms of their Schwartz kernels \cite{ioos-lu-ma-ma,bergman}, one can show that the operator $\Delta^q_pP_{\mathcal H_p}$ is a Toeplitz operator. Then the identities \eqref{e:Jqq} are easily derived from the properties of Toeplitz operators (\cite{ioos-lu-ma-ma,bergman}, in particular, \cite[Proposition 4.3]{ioos-lu-ma-ma}). See also \cite{ma-ma08a} (in particular, \cite[formula (2.31)]{ma-ma08a}) for another approach.  

Finally, by \cite[Theorem 3.1]{ma-ma08a}, the spectral density function $\rho$ is related with the coefficient $J_{1,2}$ by
\[
\rho(x_0)=J_{1,2,x_0}(0,0).
\]

\subsection{The Ma-Marinescu formula}
Here, we recall a formula for the coefficient $F_{1,2,x_0}$ proved in \cite[Subsection 2.1, formula (2.12)]{ma-ma08a}. Recall that we identify $Z\in{\mathbb R}^{2n}$ with the tangent vector $Z= \sum_{j=1}^{2n} Z_j e_j\in T_{x_0}X$. One can write 
\begin{equation}\label{e:Z=z+barz}
Z =\sum_{j=1}^nz_j\frac{\partial}{\partial z_j}+\sum_{j=1}^n\bar{z}_j\frac{\partial}{\partial \bar{z}_j},
\end{equation}
where $z_j=Z_{2j-1}+iZ_{2j}, j=1,\ldots,n$. We denote by $R^{TX}$ the curvature of $\nabla^{TX}$.  Then 
\begin{equation}\label{e:F120}
F_{1,2,x_0}(Z,Z^\prime) = [{\mathcal P}_{x_0}\mathcal F_{1,2,x_0}{\mathcal P}_{x_0}] (Z,Z^\prime),
\end{equation}
where $\mathcal F_{1,2,x_0}$ is an unbounded linear operator in $L^2(T_{x_0}X)$ given by 
\begin{multline}\label{e:F12}
\mathcal F_{1,2,x_0}= 4 \sum_{j,k=1}^n\left \langle R^{TX}_{x_0} \left(\tfrac{\partial}{\partial z_j}, \tfrac{\partial}{\partial z_k}\right) \tfrac{\partial}{\partial \bar{z}_j}, \tfrac{\partial}{\partial \bar{z}_k}\right \rangle \\
\begin{aligned}
& + \sum_{j=1}^n\left \langle (\nabla ^{X} \nabla ^{X}\mathcal J)_{(Z,Z)} \tfrac{\partial}{\partial z_j},
 \tfrac{\partial}{\partial \bar{z}_j} \right \rangle 
-\frac{1}{4} 
{\mbox{\rm tr}}_{|TX} \Big(\nabla ^{X} \nabla ^{X}|\mathcal J|\Big)_{(Z,Z)} \\
& + \frac{1}{9} |(\nabla_Z^X \mathcal J)Z|^2
+ \frac{4}{9} \sum_{j,j^\prime=1}^n \left\langle(\nabla ^{X}_Z \mathcal J)Z, \tfrac{\partial} {\partial z_j}\right\rangle b^+_j {\mathcal L}^{-1}b_{j^\prime} \left\langle(\nabla ^{X}_Z\mathcal J)Z,\tfrac{\partial} {\partial\bar{z}_{j^\prime}}\right\rangle .
\end{aligned}
\end{multline}
Here and later on, all tensors are evaluated at $x_0$, and we omit the subscript $x_0$. We identify a function $f(Z)$ on $T_{x_0}X\cong {\mathbb R}^{2n}$ with the corresponding multiplication operator in $L^2(T_{x_0}X)$.

\subsection{Preliminary considerations}
By \eqref{e:def-omega}, \eqref{e:defJ0}, and \eqref{e:defcJ0}, 
\begin{equation}\label{e:JVW}
\langle \mathcal JV,W\rangle =R^L(V,W), \quad V,W\in TX. 
\end{equation}
It follows that
\begin{equation}\label{e:JVW1}
\langle (\nabla^X_U\mathcal J)V,W\rangle =(\nabla^X_UR^L)(V,W), \quad U,V,W\in TX. 
\end{equation}
Using the fact that the 2-form $R^L$ is closed, one can show that 
\begin{equation}\label{e:JVW2}
\langle (\nabla^X_U\mathcal J)V,W\rangle+\langle (\nabla^X_V\mathcal J)W,U\rangle+\langle (\nabla^X_W\mathcal J)U,V\rangle=0, \quad U,V,W\in TX.
\end{equation}
Note that, by \eqref{e:JVW}, we easily get
\[
\langle |\mathcal J|U,V\rangle =iR^L(U,JV),\quad U,V\in TX.
\]

For any $j=1,\ldots,n$, we define a bilinear form $q_j$ on $T_{x_0}X\cong {\mathbb R}^{2n}$ by 
\[
q_j(U,V)= \left\langle(\nabla ^{X}_{U}\mathcal J)V,\tfrac{\partial} {\partial\bar{z}_j}\right\rangle, \quad U,V\in T_{x_0}X. 
\]
Since $\mathcal J$ is purely imaginary (cf. \eqref{e:JVW1}), we have
\[
\overline{q_j(U,V)} = -\left\langle(\nabla ^{X}_{U} \mathcal J)V, \tfrac{\partial} {\partial z_j}\right\rangle, \quad U,V\in T_{x_0}X.
\]
Using this notation, we can write \eqref{e:F12} as follows:
\begin{multline*}
\mathcal F_{1,2,x_0}
= 4 \sum_{j,k=1}^n \left \langle R^{TX}_{x_0} \left(\tfrac{\partial}{\partial z_j}, \tfrac{\partial}{\partial z_k}\right) \tfrac{\partial}{\partial \bar{z}_j}, 
\tfrac{\partial}{\partial \bar{z}_k}\right \rangle\\
+ Q_2(Z,Z) + Q_4(Z) - \frac{4}{9} \sum_{j,j^\prime=1}^n\overline{q_j(Z,Z)}b^+_j {\mathcal L}^{-1}b_{j^\prime} q_{j^\prime}(Z,Z),
\end{multline*}
where $Q_2$ is a bilinear form on $T_{x_0}X$ given by
\[
Q_2(U,V)=\sum_{j=1}^n \left \langle (\nabla ^{X} \nabla ^{X}\mathcal J)_{(U,V)} \tfrac{\partial}{\partial z_j}, \tfrac{\partial}{\partial \bar{z}_j} \right \rangle 
-\frac{1}{4} {\mbox{\rm tr}}_{|TX} \Big(\nabla ^{X} \nabla ^{X}|\mathcal J|\Big)_{(U,V)}, \quad U,V\in T_{x_0}X.
\]
and $Q_4$ is given by 
\[
Q_4(Z)= \frac{1}{9}  |(\nabla_Z^X \mathcal J) Z|^2=\frac{4}{9} \sum_{j=1}^n\left| q_j(Z,Z) \right|^2, \quad Z\in T_{x_0}X. 
\]

Accordingly, we have 
\begin{equation}\label{e:F120a}
J_{1,2}(0,0) =A_0+A_1+A_2-A_3,
\end{equation}
where 
\begin{align*}
A_0= & 4\sum_{j,k=1}^n \left \langle R^{TX}_{x_0} \left(\tfrac{\partial}{\partial z_j}, \tfrac{\partial}{\partial z_k}\right) \tfrac{\partial}{\partial \bar{z}_j}, 
\tfrac{\partial}{\partial \bar{z}_k}\right \rangle, \\
A_1=& \frac{1}{\mathcal P(0, 0)}[{\mathcal P} Q_2{\mathcal P}] (0,0), \quad  
A_2=\frac{1}{\mathcal P(0, 0)}[{\mathcal P} Q_4{\mathcal P}] (0,0),\\
A_3=& \frac{4}{9\mathcal P(0, 0)}\sum_{j,j^\prime=1}^n [{\mathcal P}\overline{q_j(Z,Z)} b^+_j {\mathcal L}^{-1}b_{j^\prime} q_{j^\prime}(Z,Z){\mathcal P}] (0,0)
\end{align*}

By \eqref{e:Z=z+barz}, we can write 
\begin{equation}\label{e:qj}
q_j(Z,Z)=\sum_{k,\ell=1}^n (q_{j,k\ell}z_kz_\ell + q_{j,k\bar \ell}z_k\bar z_\ell+ q_{j,\bar k\bar \ell}\bar z_k\bar z_\ell),
\end{equation}
where
\begin{align*}
q_{j,k\ell}=& q_j\left(\tfrac{\partial}{\partial z_k}, \tfrac{\partial}{\partial z_\ell}\right)=\left\langle(\nabla ^{X}_{\tfrac{\partial}{\partial z_k}}\mathcal J)\tfrac{\partial}{\partial z_\ell},\tfrac{\partial} {\partial\bar{z}_j}\right\rangle,\\
q_{j,k\bar \ell}=& q_j\left(\tfrac{\partial}{\partial z_k}, \tfrac{\partial}{\partial \bar z_\ell}\right)+q_j\left(\tfrac{\partial}{\partial \bar z_\ell}, \tfrac{\partial}{\partial z_k}\right) = \left\langle(\nabla ^{X}_{\tfrac{\partial}{\partial z_k}}\mathcal J)\tfrac{\partial}{\partial \bar z_\ell},\tfrac{\partial} {\partial\bar{z}_j}\right\rangle+\left\langle(\nabla ^{X}_{\tfrac{\partial}{\partial \bar z_\ell}}\mathcal J)\tfrac{\partial}{\partial z_k},\tfrac{\partial} {\partial\bar{z}_j}\right\rangle, \\
q_{j,\bar k\bar \ell}=& q_j\left(\tfrac{\partial}{\partial \bar z_k}, \tfrac{\partial}{\partial \bar z_\ell}\right)=\left\langle(\nabla ^{X}_{\tfrac{\partial}{\partial \bar z_k}}\mathcal J)\tfrac{\partial}{\partial \bar z_\ell},\tfrac{\partial} {\partial\bar{z}_j}\right\rangle.
\end{align*}

Using the fact that $\nabla ^{X}_{\partial/\partial \bar z_k}\mathcal J$ is skew-adjoint and \eqref{e:JVW2}, one can easily see that
\begin{equation}\label{e:q-ellkj1}
q_{\ell,\bar k\bar j}=-q_{j,\bar k\bar \ell},\quad q_{j,\bar k\bar j}=0, \quad q_{j,\bar k\bar \ell}+q_{\ell,\bar j\bar k}+q_{k,\bar \ell\bar j}=0. 
\end{equation}
From the equality $|q_{j,\bar k\bar \ell}-q_{j,\bar \ell\bar k}|^2=|q_{k,\bar j\bar \ell}|^2$, we infer
\begin{equation}\label{e:q-ellkj2}
\overline{q_{j,\bar k\bar \ell}}q_{j,\bar \ell\bar k}+q_{j,\bar k\bar \ell}\overline{q_{j,\bar \ell\bar k}}=|q_{j,\bar k\bar \ell}|^2+|q_{j,\bar \ell\bar k}|^2-|q_{k,\bar j\bar \ell}|^2.
\end{equation}
and 
\begin{equation}\label{e:q-ellkj3}
|q_{j,\bar k\bar \ell}+q_{j,\bar \ell\bar k}|^2=2|q_{j,\bar k\bar \ell}|^2+2|q_{j,\bar \ell\bar k}|^2-|q_{k,\bar j\bar \ell}|^2. 
\end{equation}

\subsection{Eigenfunctions and matrix  coefficients of the model operator}
Here, we recall some necessary information on eigenvalues and  eigenfunctions of the model operator $\mathcal L$ (cf. \cite[Section 1.4]{ma-ma08a}) and make some computations.  

For any polynomial  $g(z,\bar{z})$ on $z$ and $\bar{z}$
\begin{gather}
[b_i,b^{+}_j]=b_i b^{+}_j-b^{+}_j b_i =-2a_i \delta_{i\,j},\quad 
[b_i,b_j]=[b^{+}_i,b^{+}_j]=0\, ,\label{e:com1}\\
[g(z,\bar{z}),b_j]=  2 \tfrac{\partial}{\partial z_j}g(z,\bar{z}), 
\quad  [g(z,\bar{z}),b_j^+]
= - 2\tfrac{\partial}{\partial \bar{z}_j}g(z,\bar{z})\,. \label{com-bg}
\end{gather}

By \cite[formula (1.98)]{ma-ma08a}, 
\begin{equation}\label{bP}
(b^{+}_j\mathcal P)(Z,Z^\prime)=0, \quad (b_j\mathcal P)(Z,Z^\prime)=a_j(\bar z_j-\bar z^\prime_j)\mathcal P(Z,Z^\prime).
\end{equation}
For $|\alpha|>0$ and for any polynomial $g(z,\bar z)$ we have
\[
\mathcal P b^\alpha g(z,\bar z)\mathcal P=0.
\]

We recall \cite[Theorem 1.15]{ma-ma08a} that the spectrum of ${\mathcal L}$ on $L^2({\mathbb R}^{2n})$ is given by
\begin{equation}\label{bk2.68}
\sigma({\mathcal L})=
\left(2\sum_{i=1}^n\alpha_i a_i\,:\, 
 \alpha =(\alpha_1,\cdots,\alpha_n)\in {\mathbb Z}_+^n\right)
\end{equation}
and an orthogonal base of the eigenspace of $2\sum_{i=1}^n\alpha_i a_i$
is given by
\begin{equation}\label{bk2.69}
\Phi_{\alpha,\beta}=b^{\alpha}\left(z^{\beta}\exp\left({-\frac{1}{4}\sum_{i=1}^n
a_i|z_i|^2}\right)\right)\,,\quad \beta\in{\mathbb Z}_+^n.
\end{equation}
So we have
\begin{equation}\label{eigenL0}
{\mathcal L} \Phi_{\alpha,\beta}=\left(2\sum_{i=1}^n\alpha_i a_i\right)\Phi_{\alpha,\beta}, \quad \alpha, \beta\in{\mathbb Z}_+^n.
\end{equation}
In particular, an orthonormal base of $\operatorname{ker} ({\mathcal L})$ is formed by the functions
\begin{equation}\label{bk2.70}
\varphi_\beta(z)=\left(\frac{a ^\beta}{(2\pi)^n 2 ^{|\beta|} \beta!}
\prod_{i=1}^n a_i\right)^{1/2}z^\beta
\exp\Big (-\frac{1}{4} \sum_{j=1}^n a_j |z_j|^2\Big )\,,\quad \beta\in{\mathbb Z}_+^n.
\end{equation}
and for $|\alpha|>0$
\begin{equation}\label{L0-1}
{\mathcal L}^{-1} \Phi_{\alpha,\beta}=\frac{1}{2\sum_{i=1}^n\alpha_i a_i}\Phi_{\alpha,\beta}.
\end{equation}
Observe that for any $\beta\in{\mathbb Z}_+^n$
\[ 
z^\beta\mathcal P(Z, 0)=\left(\frac{2 ^{|\beta|} \beta!}{(2\pi)^n a ^\beta}\prod_{j=1}^na_j\right)^{1/2} \varphi_\beta(z).
\]
In particular, we have
\[
\mathcal P(0, 0)=\frac{1}{(2\pi)^n}\prod_{j=1}^na_j. 
\]
We use this relation to simplify formulas. 

Thus, for the $L^2$-norms, we have
\begin{equation}\label{e:L2-norm}
\|z^\beta\mathcal P(Z, 0)\|^2_{L^2}=\frac{2 ^{|\beta|} \beta!}{a ^\beta}\mathcal P(0, 0), \quad \beta\in{\mathbb Z}_+^n.
\end{equation}
It follows that
\begin{equation}\label{e:L2-norm2}
\|z^\beta\bar z^{\gamma}\mathcal P(Z, 0)\|^2_{L^2}=\frac{2 ^{|\beta+\gamma|} (\beta+\gamma)!}{a^{\beta+\gamma}}\mathcal P(0, 0), \quad \beta,\gamma\in{\mathbb Z}_+^n.
\end{equation}
Using \eqref{e:com1}, \eqref{com-bg} and \eqref{bP}, one can easily derive the following formulas for any $j,k=1,\ldots, n$ and $\beta\in{\mathbb Z}_+^n$:
\begin{align}\label{e:L2-norm3}
\|b_j(z^\beta\mathcal P(Z, 0))\|^2_{L^2} & = 2a_j\frac{2 ^{|\beta|} \beta!}{a ^\beta}\mathcal P(0, 0),\\
\label{e:L2-norm4}
\|b_jb_k(z^\beta\mathcal P(Z, 0))\|^2_{L^2} & =4(1+\delta_{jk})a_ja_k \frac{2 ^{|\beta|} \beta!}{a ^\beta}\mathcal P(0, 0).
\end{align}

\section{Proof of Theorem \ref{t:main}}
In this section, we compute the terms of \eqref{e:F120a} and complete the proof of Theorem \ref{t:main}.

\subsection{Computation of $A_0$} Using \eqref{e:eigencJ} and the fact that $\mathcal J$ is skew-adjoint, we proceed as follows:
\begin{multline*}
\left \langle R^{TX}_{x_0} \left(\tfrac{\partial}{\partial z_j}, \tfrac{\partial}{\partial z_k}\right) \tfrac{\partial}{\partial \bar{z}_j},  \tfrac{\partial}{\partial \bar{z}_k}\right \rangle\\
=-\frac{1}{a_k} \left \langle R^{TX}_{x_0} \left(\tfrac{\partial}{\partial z_j}, \tfrac{\partial}{\partial z_k}\right) \tfrac{\partial}{\partial \bar{z}_j}, 
\mathcal J \tfrac{\partial}{\partial \bar{z}_k}\right \rangle
= \frac{1}{a_k} \left \langle \mathcal J R^{TX}_{x_0} \left(\tfrac{\partial}{\partial z_j}, \tfrac{\partial}{\partial z_k}\right) \tfrac{\partial}{\partial \bar{z}_j}, 
\tfrac{\partial}{\partial \bar{z}_k}\right \rangle\\
= -\frac{a_j}{a_k} \left \langle R^{TX}_{x_0} \left(\tfrac{\partial}{\partial z_j}, \tfrac{\partial}{\partial z_k}\right) \tfrac{\partial}{\partial \bar{z}_j},  \tfrac{\partial}{\partial \bar{z}_k}\right \rangle  - \frac{1}{a_k} \left \langle \left[R^{TX}_{x_0} \left(\tfrac{\partial}{\partial z_j}, \tfrac{\partial}{\partial z_k}\right),\mathcal J\right] \tfrac{\partial}{\partial \bar{z}_j}, 
\tfrac{\partial}{\partial \bar{z}_k}\right \rangle.
\end{multline*}
We infer that
\begin{equation}\label{e:RR}
\left \langle R^{TX}_{x_0} \left(\tfrac{\partial}{\partial z_j}, \tfrac{\partial}{\partial z_k}\right) \tfrac{\partial}{\partial \bar{z}_j},  \tfrac{\partial}{\partial \bar{z}_k}\right \rangle 
 =-\frac{1}{a_j+a_k} \left \langle \left[R^{TX}_{x_0} \left(\tfrac{\partial}{\partial z_j}, \tfrac{\partial}{\partial z_k}\right), \mathcal J\right] \tfrac{\partial}{\partial \bar{z}_j}, 
\tfrac{\partial}{\partial \bar{z}_k}\right \rangle.
\end{equation}
Recall the following identity (\cite[formula (2.2)]{ma-ma08a}): for any $A\in\operatorname{End}(TX)$
\begin{equation}\label{e:comm-nabla}
(\nabla ^{X} \nabla ^{X}A)_{(U,V)}-(\nabla ^{X} \nabla ^{X}A)_{(V,U)}=\left[R^{TX}(U,V),A\right].
\end{equation}
By \eqref{e:RR} and \eqref{e:comm-nabla}, we get
\begin{multline}
A_0 = -\sum_{j,k=1}^n \frac{4}{a_j+a_k} \Big \langle \Big( (\nabla ^{X} \nabla ^{X}\mathcal J)_{(\tfrac{\partial}{\partial z_j},\tfrac{\partial}{\partial z_k})}-(\nabla ^{X} \nabla ^{X}\mathcal J)_{(\tfrac{\partial}{\partial z_k},\tfrac{\partial}{\partial z_j})}\Big) \tfrac{\partial}{\partial \bar{z}_j}, 
\tfrac{\partial}{\partial \bar{z}_k} \Big \rangle\\ 
= -\sum_{j,k=1}^n\frac{8}{a_j+a_k} \left \langle (\nabla ^{X} \nabla ^{X}\mathcal J)_{(\tfrac{\partial}{\partial z_j},\tfrac{\partial}{\partial z_k})} \tfrac{\partial}{\partial \bar{z}_j}, 
\tfrac{\partial}{\partial \bar{z}_k} \right \rangle.\label{e:A0}
\end{multline}

\subsection{Computation of $A_1$} 
Let $q$ be a quadratic function on ${\mathbb R}^{2n}$ of the form $q(Z)=Q(Z,Z)$ with some bilinear (not necessarily, symmetric) form $Q$. As above, we can write 
\[
q(Z)=\sum_{k,\ell=1}^n (q_{k\ell}z_kz_\ell + q_{k\bar \ell}z_k\bar z_\ell+ q_{\bar k\bar \ell}\bar z_k\bar z_\ell),
\]
where
\[
q_{k\ell}=Q\left(\tfrac{\partial}{\partial z_k}, \tfrac{\partial}{\partial z_\ell}\right),\quad
q_{\bar k\bar \ell}=Q\left(\tfrac{\partial}{\partial \bar z_k}, \tfrac{\partial}{\partial \bar z_\ell}\right), \quad
q_{k\bar \ell}=Q\left(\tfrac{\partial}{\partial z_k}, \tfrac{\partial}{\partial \bar z_\ell}\right)+Q\left(\tfrac{\partial}{\partial \bar z_\ell}, \tfrac{\partial}{\partial z_k}\right),
\]
Using \eqref{bP}, \eqref{com-bg} and \eqref{eigenL0}, it is easy to see that 
\[
\mathcal P z_k z_\ell\mathcal P(0,0)=\mathcal P \bar z_k \bar z_\ell\mathcal P(Z,0)=0, \quad \mathcal P z_k \bar z_\ell\mathcal P(Z,0)=\frac{2}{a_\ell} \delta_{k\ell} \mathcal P(Z,0).
\]
Therefore, we get
\begin{equation}\label{e:PqP}
[\mathcal P q(Z) \mathcal P](0,0)=\sum_{j=1}^n \frac{2}{a_j} q_{j\bar j}\mathcal P(0, 0).
\end{equation}

By \eqref{e:comm-nabla}, \eqref{e:eigencJ}, and the fact that $\mathcal J$ is skew-adjoint, we have
\begin{multline*}
\left \langle (\nabla ^{X} \nabla ^{X}\mathcal J)_{(\tfrac{\partial}{\partial z_j},\tfrac{\partial}{\partial \bar z_j})} \tfrac{\partial}{\partial z_k}, \tfrac{\partial}{\partial \bar{z}_k} \right \rangle - \left \langle (\nabla ^{X} \nabla ^{X}\mathcal J)_{(\tfrac{\partial}{\partial \bar z_j},\tfrac{\partial}{\partial z_j})} \tfrac{\partial}{\partial z_k},
 \tfrac{\partial}{\partial \bar{z}_k} \right \rangle\\ =\left \langle \left[R^{TX}\left(\tfrac{\partial}{\partial z_j},\tfrac{\partial}{\partial \bar z_j}\right),\mathcal J\right] \tfrac{\partial}{\partial z_k}, \tfrac{\partial}{\partial \bar{z}_k} \right \rangle \\ =\left \langle R^{TX}\left(\tfrac{\partial}{\partial z_j},\tfrac{\partial}{\partial \bar z_j}\right)\mathcal J\tfrac{\partial}{\partial z_k}, \tfrac{\partial}{\partial \bar{z}_k} \right \rangle+\left \langle  R^{TX}\left(\tfrac{\partial}{\partial z_j},\tfrac{\partial}{\partial \bar z_j}\right) \tfrac{\partial}{\partial z_k}, \mathcal J \tfrac{\partial}{\partial \bar{z}_k} \right \rangle=0.
\end{multline*}
Similarly, by \eqref{e:comm-nabla}, \eqref{e:eigen|J|}, and the fact that $|\mathcal J|$ is self-adjoint, one can show that
\[
\left \langle (\nabla ^{X} \nabla ^{X}|\mathcal J|)_{(\tfrac{\partial}{\partial z_j},\tfrac{\partial}{\partial \bar z_j})} \tfrac{\partial}{\partial z_k}, \tfrac{\partial}{\partial \bar{z}_k} \right \rangle - \left \langle (\nabla ^{X} \nabla ^{X}|\mathcal J|)_{(\tfrac{\partial}{\partial \bar z_j},\tfrac{\partial}{\partial z_j})} \tfrac{\partial}{\partial z_k}, \tfrac{\partial}{\partial \bar{z}_k} \right \rangle=0.
\]
Observe that for any $A\in \operatorname{End}(T_{x_0}X)$
\[
{\mbox{\rm tr}}_{|TX} A=4\sum_{k=1}^n \left \langle A \tfrac{\partial}{\partial z_k}, \tfrac{\partial}{\partial \bar{z}_k} \right \rangle.
\]
It follows that 
\[
Q_2\left(\tfrac{\partial}{\partial z_j},\tfrac{\partial}{\partial \bar z_j}\right)=Q_2\left(\tfrac{\partial}{\partial \bar z_j},\tfrac{\partial}{\partial z_j}\right), 
\]
and, by \eqref{e:PqP}, we get 
\begin{equation}\label{e:A1}
A_1 = \sum_{j=1}^n \frac{1}{a_j} \left( {\mbox{\rm tr}}_{|TX} (\nabla ^{X} \nabla ^{X}\mathcal J)_{(\tfrac{\partial}{\partial z_j},\tfrac{\partial}{\partial \bar z_j})}
 - {\mbox{\rm tr}}_{|TX} \Big(\nabla ^{X} \nabla ^{X}|\mathcal J|\Big)_{(\tfrac{\partial}{\partial z_j},\tfrac{\partial}{\partial \bar z_j})}\right).
\end{equation}

\subsection{Computation of $A_2$}
One can write
\begin{multline*}
A_2=\frac 49 \frac{1}{\mathcal P(0, 0)}\sum_{j=1}^n[\mathcal P |q_j(Z,Z)|^2 \mathcal P](0,0)\\ =\frac 49  \frac{1}{\mathcal P(0, 0)}\sum_{j=1}^n \int \mathcal P(0,Z) |q_j(Z,Z)|^2 \mathcal P(Z,0)dZ=\frac 49  \frac{1}{\mathcal P(0, 0)}\sum_{j=1}^n\|q_j(Z,Z)\mathcal P(Z,0)\|_{L^2}^2. 
\end{multline*}
By \eqref{e:qj}, \eqref{bP} and \eqref{com-bg}, for any $j=1,2,\ldots,n$, we have 
\begin{align*}
q_j(Z,Z)\mathcal P(Z,0) = & \sum_{k,\ell=1}^n (q_{j,k\ell}z_kz_\ell \mathcal P(Z,0) + q_{j,k\bar \ell}z_k\bar z_\ell \mathcal P(Z,0) + q_{j,\bar k\bar \ell}\bar z_k\bar z_\ell\mathcal P(Z,0)) \\
= & \sum_{k,\ell=1}^n q_{j,k\ell}z_kz_\ell \mathcal P(Z,0) +\sum_{k,\ell=1}^n \frac{1}{a_\ell }q_{j,k\bar \ell}b_\ell z_k \mathcal P(Z,0)\\ & + \sum_{k,\ell=1}^n \frac{2}{a_\ell}q_{j,k\bar \ell}\delta_{k\ell} \mathcal P(Z,0) +\sum_{k,\ell=1}^n q_{j,\bar k\bar \ell}\frac{1}{a_ka_\ell} b_k b_\ell  \mathcal P(Z,0).
\end{align*}
Since $\Phi_{\alpha,\beta}$ are orthogonal, we get
\begin{align*}
\|q_j(Z,Z)\mathcal P(Z,0)\|_{L^2}^2 = & \sum_{k=1}^n |q_{j,kk}|^2\|z^2_k \mathcal P(Z,0)\|_{L^2}^2+\sum_{k<\ell} |q_{j,k\ell}+q_{j,\ell k}|^2\|z_kz_\ell \mathcal P(Z,0)\|_{L^2}^2\\ & + \sum_{k,\ell=1}^n \frac{1}{a^2_\ell } \left|q_{j,k\bar \ell}\right|^2 \| b_\ell z_k \mathcal P(Z,0)\|_{L^2}^2 + \sum_{k,\ell=1}^n \frac{4}{a^2_\ell} \left|q_{j,k\bar \ell}\right|^2\delta_{k\ell} \|\mathcal P(Z,0)\|_{L^2}^2\\ & + \sum_{k=1}^n|q_{j,\bar k\bar k}|^2\frac{1}{a^4_k} \|b^2_k\mathcal P(Z,0)\|_{L^2}^2+ \sum_{k<\ell}^n|q_{j,\bar k\bar \ell}+q_{j,\bar \ell \bar k}|^2\frac{1}{a^2_ka^2_\ell} \|b_kb_\ell \mathcal P(Z,0)\|_{L^2}^2.
\end{align*}
Using \eqref{e:L2-norm},  \eqref{e:L2-norm2}, and \eqref{e:L2-norm3}, we get
\[
\|q_j(Z,Z)\mathcal P(Z,0)\|_{L^2}^2 = \sum_{k,\ell=1}^n \frac{1}{a_ka_\ell}\left(2|q_{j,k\ell}+q_{j,\ell k}|^2+ 4(1+\delta_{k\ell}) \left|q_{j,k\bar \ell}\right|^2  + 2|q_{j.\bar k\bar \ell}+q_{j.\bar \ell\bar k}|^2 \right)\mathcal P(0, 0),
\]
and
\begin{equation} \label{e:A2}
A_2 = \sum_{j,k,\ell=1}^n \frac{8}{9a_ka_\ell}\left(|q_{j,k\ell}+q_{j,\ell k}|^2+ 2(1+\delta_{k\ell}) \left|q_{j,k\bar \ell}\right|^2  + |q_{j.\bar k\bar \ell}+q_{j.\bar \ell\bar k}|^2 \right).
\end{equation}

\subsection{Computation of $A_3$} Observe that 
\[
\mathcal P \bar q_j b^+_j {\mathcal L}^{-1}b_{j^\prime} q_{j^\prime}\mathcal P(0,0)=\int [\mathcal P \bar q_j b^+_j](0,Z) [{\mathcal L}^{-1}b_{j^\prime} q_{j^\prime}\mathcal P](Z,0)dZ 
\]
where $[\mathcal P \bar q_j b^+_j](Z,Z^\prime)$ and $[{\mathcal L}^{-1}b_{j^\prime} q_{j^\prime}\mathcal P](Z,Z^\prime)$ are the Schwartz kernels of the operators $\mathcal P \bar q_j b^+_j$ and ${\mathcal L}^{-1}b_{j^\prime} q_{j^\prime}\mathcal P$ with respect to $dZ$. The operator $\mathcal P \bar q_j b^+_j$ is the adjoint of the operator $b_j q_j \mathcal P$ in $L^2({\mathbb R}^{2n})$. Therefore, $
[\mathcal P \bar q_j b^+_j](Z,Z^\prime)=\overline{[b_j q_j \mathcal P](Z^\prime,Z)}$. It follows that
\begin{equation}\label{e:A3L2}
\mathcal P \bar q_j b^+_j {\mathcal L}^{-1}b_{j^\prime} q_{j^\prime}\mathcal P(0,0)=\left( {\mathcal L}^{-1}b_{j^\prime} q_{j^\prime}\mathcal P(Z,0), b_j q_j\mathcal P(Z,0) \right)_{L^2},
\end{equation}
where $\left(\cdot,\cdot\right)_{L^2}$ denotes the inner product in $L^2({\mathbb R}^{2n})$. Let us compute the inner product on the right-hand side of \eqref{e:A3L2}. 
By \eqref{e:qj}, \eqref{bP}, and \eqref{com-bg}, we get 
\begin{equation}\label{e:bjqj}
\begin{aligned}
b_j q_j\mathcal P(Z,0) = &  \sum_{k,\ell=1}^n q_{j,k\ell}b_j z_kz_\ell \mathcal P(Z,0) + \sum_{k,\ell=1}^n q_{j,k\bar \ell}\frac{1}{a_\ell} b_j b_\ell z_k\mathcal P(Z,0)\\  & + \sum_{k=1}^n q_{j,k\bar k}\tfrac{2}{a_k}b_j\mathcal P(Z,0) + \sum_{k,\ell=1}^n q_{j,\bar k\bar \ell}\frac{1}{a_ka_\ell} b_j b_k b_\ell  \mathcal P(Z,0)\\ = & u_1(Z)+u_2(Z)+u_3(Z)+u_4(Z).
\end{aligned}
\end{equation}
Now we use \eqref{L0-1}: 
\begin{multline}\label{e:L0-1bjqj}
{\mathcal L}^{-1} b_{j^\prime} q_{j^\prime}\mathcal P(Z,0)\\
\begin{aligned}
= & \sum_{k^\prime,\ell^\prime=1}^n q_{{j^\prime},k^\prime\ell^\prime}\frac{1}{2a_{j^\prime}} b_{j^\prime} z_{k^\prime}z_{\ell^\prime} \mathcal P(Z,0) +  \sum_{k^\prime,\ell^\prime=1}^n q_{{j^\prime},k^\prime\bar \ell^\prime}\frac{1}{a_{\ell^\prime}}\frac{1}{2(a_{j^\prime}+a_{\ell^\prime})} b_{j^\prime} b_{\ell^\prime} z_{k^\prime}\mathcal P(Z,0) \\ 
&   +  \sum_{k^\prime=1}^n q_{{j^\prime},k^\prime\bar k^\prime}\frac{1}{a_{k^\prime} a_{j^\prime}}b_{j^\prime}\mathcal P(Z,0)+  \sum_{k^\prime,\ell^\prime=1}^n q_{{j^\prime},\bar k^\prime\bar \ell^\prime}\frac{1}{a_{k^\prime}a_{\ell^\prime}} \frac{1}{2(a_{j^\prime}+a_{k^\prime}+a_{\ell^\prime})} b_{j^\prime} b_{k^\prime} b_{\ell^\prime}  \mathcal P(Z,0) \\
= & v_1(Z)+v_2(Z)+v_3(Z)+v_4(Z).
\end{aligned} 
\end{multline}
We use \eqref{e:bjqj} and \eqref{e:L0-1bjqj} to compute the right-hand side of \eqref{e:A3L2}. Since $\Phi_{\alpha,\beta}$ are orthogonal, the terms $u_j, j=1,2,3,4,$ on the right-hand side of \eqref{e:bjqj} (accordingly, the terms $v_j, j=1,2,3,4,$ on the right-hand side of \eqref{e:L0-1bjqj}) are mutually orthogonal. Therefore,  the right-hand side of \eqref{e:A3L2} is the sum of four terms, each term is the inner product of the corresponding terms on the right-hand side of \eqref{e:bjqj} and \eqref{e:L0-1bjqj}:
\begin{equation} \label{e:A3}
A_3=\frac{4}{9\mathcal P(0, 0)}(I_1+I_2+I_3+I_4), \quad I_j=(u_j,v_j)_{L^2}, j=1,2,3,4.
\end{equation}
We compute these inner products. For the first three of them we get
\begin{align*}
I_1 = & \sum_{j,k=1}^n|q_{j,kk}|^2\frac{1}{2a_j}\| b_j z^2_k \mathcal P(Z,0)\|^2_{L^2} +\sum_{j=1}^n \sum_{k<\ell}|q_{j,k\ell}+q_{j,\ell k}|^2\frac{1}{2a_j}\| b_j z_kz_\ell \mathcal P(Z,0)\|^2_{L^2}\\= & \sum_{j,k,\ell=1}^n |q_{j,k\ell}+q_{j,\ell k}|^2\frac{2}{a_ka_\ell}\mathcal P(0, 0);\\
I_2= & \sum_{j,k=1}^n |q_{j,k\bar j}|^2\frac{1}{4a^3_j}\| b^2_j z_k\mathcal P(Z,0)\|^2_{L^2} +\sum_{k=1}^n \sum_{j\neq\ell} \left|q_{j,k\bar \ell}\frac{1}{a_\ell}+q_{\ell,k\bar j}\frac{1}{a_j}\right|^2\frac{1}{4(a_j+a_\ell)}\| b_j b_\ell z_k\mathcal P(Z,0)\|^2_{L^2}\\  = & \sum_{j,k=1}^n |q_{j,k\bar j}|^2\frac{4}{a_ka_j}\mathcal P(0,0) + \sum_{k=1}^n \sum_{j\neq\ell} \left|q_{j,k\bar \ell}\frac{1}{a_\ell}+q_{\ell,k\bar j}\frac{1}{a_j}\right|^2\frac{2a_ja_\ell}{a_k(a_j+a_\ell)}\mathcal P(0, 0)\\  
= & \sum_{j,k,\ell=1}^n \left(|q_{j,k\bar \ell}|^2\frac{4a_j}{a_ka_\ell(a_j+a_\ell)}+(\overline{q_{j,k\bar \ell}}q_{\ell,k\bar j}+q_{j,k\bar \ell}\overline{q_{\ell,k\bar j}})\frac{2}{a_k(a_j+a_\ell)}\right)\mathcal P(0, 0);\\
I_3= & \sum_{k=1}^n |q_{j,k\bar k}|^2\frac{2}{a^2_k a_j}\|b_j\mathcal P(Z,0)\|^2_{L^2} =\sum_{k=1}^n |q_{j,k\bar k}|^2\frac{4}{a^2_k}\mathcal P(0, 0).
\end{align*}
Finally, for $I_4$, we can write
\begin{align*}
I_4= &   \sum_{j\neq k} \left|q_{j,\bar j\bar k}\frac{1}{a_ja_k}+q_{j,\bar k\bar j}\frac{1}{a_ka_j}+q_{k,\bar j\bar j}\frac{1}{a^2_j}\right|^2 \frac{1}{2(2a_j+a_k)} \|b^2_j b_k  \mathcal P(Z,0)\|^2 \\
& + \sum_{j<k<\ell} \left| \sum_{j^\prime,k^\prime,\ell^\prime}  q_{j^\prime,\bar k^\prime\bar \ell^\prime}\frac{1}{a_{k^\prime}a_{\ell^\prime}}\right|^2 \frac{1}{2(a_j+a_k+a_\ell)} \|b_j b_k b_\ell  \mathcal P(Z,0)\|^2,
\end{align*}
where the last sum is taken over all permutations $(j^\prime,k^\prime,\ell^\prime)$ of $j,k,\ell$. 
Taking into account \eqref{e:q-ellkj1}, we get
\begin{align*}
I_4 =  &   \sum_{j,k=1}^n \left|q_{j,\bar j\bar k}\frac{1}{a_ja_k}+q_{j,\bar k\bar j}\frac{1}{a_ka_j}+q_{k,\bar j\bar j}\frac{1}{a^2_j}\right|^2 \frac{8 a_j^2a_k}{2a_j+a_k} \mathcal P(0,0) \\ & + \sum_{j<k<\ell} \left| \sum_{j^\prime,k^\prime,\ell^\prime}  q_{j^\prime,\bar k^\prime\bar \ell^\prime}\frac{1}{a_{k^\prime}a_{\ell^\prime}}\right|^2 \frac{4a_ja_ka_\ell}{a_j+a_k+a_\ell} \mathcal P(0,0)\\
=  &   \sum_{j,k=1}^n |q_{k,\bar j\bar j}|^2\frac{8(a_j-a_k)^2}{a_j^2a_k(2a_j+a_k)} \mathcal P(0,0) \\ & + \sum_{j \neq k \neq \ell\neq j} \left|q_{j,\bar k\bar \ell}(a_j-a_\ell)+q_{\ell,\bar j\bar k}(a_\ell-a_k)+q_{k,\bar \ell\bar j}(a_k-a_j)\right|^2 \frac{2}{3a_ja_ka_\ell(a_j+a_k+a_\ell)} \mathcal P(0,0)\\
=  &  \sum_{j,k,\ell=1}^n |q_{j,\bar k\bar \ell}|^2\frac{2(a_j-a_\ell)^2}{a_ja_ka_\ell(a_j+a_k+a_\ell)} \mathcal P(0,0)\\ & - \sum_{j,k,\ell=1}^n(\overline{q_{j,\bar k\bar \ell}}q_{j,\bar \ell\bar k}+q_{j,\bar k\bar \ell}\overline{q_{j,\bar \ell\bar k}}) \frac{2(a_j-a_\ell)(a_k-a_j)}{a_ja_ka_\ell(a_j+a_k+a_\ell)} \mathcal P(0,0).
\end{align*}

Now we compute $A_2-A_3$. It is easy to see that the terms, containing $q_{j,k\ell}$. cancel:
\[
\sum_{j,k,\ell=1}^n \frac{8}{9a_ka_\ell}|q_{j,k\ell}+q_{j,\ell k}|^2-\frac{4}{9\mathcal P(0, 0)}I_1=0.
\]
Therefore, we get
\begin{equation} \label{e:A2-A3}
A_2-A_3=J_1+J_2, 
\end{equation}
where $J_1$ contains the terms with $q_{j,k\bar\ell}$ and $J_2$ with $q_{j,\bar k\bar\ell}$. For $J_1$, we have 
\[
J_1=\sum_{j,k,\ell=1}^n \frac{16}{9a_ka_\ell}(1+\delta_{k\ell}) \left|q_{j,k\bar \ell}\right|^2-\frac{4}{9\mathcal P(0, 0)}(I_2+I_3).
\]
Using the above expressions for $I_2$ and $I_3$, we compute:
\begin{multline*}
J_1=  \sum_{j,k,\ell=1}^n \Big( \frac{16}{9a_ka_\ell}\left|q_{j,k\bar \ell}\right|^2-|q_{j,k\bar \ell}|^2\frac{16a_j}{9a_ka_\ell(a_j+a_\ell)} -(\overline{q_{j,k\bar \ell}}q_{\ell,k\bar j}+q_{j,k\bar \ell}\overline{q_{\ell,k\bar j}})\frac{8}{9a_k(a_j+a_\ell)}\Big) \\ 
= \sum_{j,k,\ell=1}^n \frac{8}{9a_k(a_j+a_\ell)}\left|q_{j,k\bar \ell}-q_{\ell,k\bar j}\right|^2.
\end{multline*}
Since, by \eqref{e:JVW2}, we have
\begin{multline*}
q_{j,k\bar \ell}-q_{\ell,k\bar j}=\left\langle(\nabla ^{X}_{\tfrac{\partial}{\partial \bar z_\ell}}\mathcal J)\tfrac{\partial}{\partial z_k},\tfrac{\partial} {\partial\bar{z}_j}\right\rangle-\left\langle(\nabla ^{X}_{\tfrac{\partial}{\partial \bar z_j}}\mathcal J)\tfrac{\partial}{\partial z_k},\tfrac{\partial} {\partial\bar{z}_\ell}\right\rangle \\ =\left\langle(\nabla ^{X}_{\tfrac{\partial}{\partial \bar z_\ell}}\mathcal J)\tfrac{\partial}{\partial z_k},\tfrac{\partial} {\partial\bar{z}_j}\right\rangle+\left\langle(\nabla ^{X}_{\tfrac{\partial}{\partial \bar z_j}}\mathcal J)\tfrac{\partial}{\partial \bar{z}_\ell},\tfrac{\partial} {\partial z_k}\right\rangle=-\left\langle  (\nabla ^{X}_{\tfrac{\partial}{\partial z_k}}\mathcal J )\tfrac{\partial}{\partial \bar z_j},\tfrac{\partial} {\partial \bar z_\ell}\right\rangle,
\end{multline*}
we conclude that
\begin{equation} \label{e:J1}
J_1=\sum_{j,k,\ell=1}^n \frac{8}{9a_k(a_j+a_\ell)}\left|\left\langle  (\nabla ^{X}_{\tfrac{\partial}{\partial z_k}}\mathcal J )\tfrac{\partial}{\partial \bar z_j},\tfrac{\partial} {\partial \bar z_\ell}\right\rangle\right|^2.
\end{equation}

For $J_2$, we have
\[
J_2=\sum_{j,k,\ell=1}^n \frac{8}{9a_ka_\ell}|q_{j.\bar k\bar \ell}+q_{j.\bar \ell\bar k}|^2-\frac{4}{9\mathcal P(0, 0)}I_4.
\]
Using the above expression for $I_4$, \eqref{e:q-ellkj2} and \eqref{e:q-ellkj3}, after routine computations, we get 
\begin{equation} \label{e:J2}
J_2= \sum_{j, k,\ell=1}^n |q_{j.\bar k\bar \ell}|^2\frac{8}{a_k(a_j+a_k+a_\ell)}.
\end{equation}
Combining \eqref{e:A0}, \eqref{e:A1}, \eqref{e:A2-A3}, \eqref{e:J1} and \eqref{e:J2}, we obtain \eqref{e:rho}. 

\subsection{Some particular cases}
One can write the main formula \eqref{e:rho} in terms of the polar decomposition  $\mathcal J=-iJ |\mathcal J|$ of $\mathcal J$. 
By routine computations, which we omit, we get:
\begin{multline*} 
\rho(x_0) = \sum_{j,k=1}^n \frac{2a_k(a_k-a_j)}{a_j(a_j+a_k)} \left \langle 
\left ( (\nabla ^{X}_{\tfrac{\partial}{\partial \bar z_j}}J) (\nabla ^{X}_{\tfrac{\partial}{\partial z_j}}J)+ (\nabla ^{X}_{\tfrac{\partial}{\partial z_j}}J) (\nabla ^{X}_{\tfrac{\partial}{\partial \bar z_j}}J)\right)\tfrac{\partial}{\partial z_k},  \tfrac{\partial}{\partial \bar z_k}  
\right \rangle \\ 
\begin{aligned}
 & +2\sum_{j,k=1}^n \left \langle \left( (\nabla ^{X}_{\tfrac{\partial}{\partial z_j}}J) (\nabla ^{X}_{\tfrac{\partial}{\partial \bar z_k}}J)+ (\nabla ^{X}_{\tfrac{\partial}{\partial \bar z_k}}J) (\nabla ^{X}_{\tfrac{\partial}{\partial z_j}}J)\right) \tfrac{\partial}{\partial \bar z_j},  \tfrac{\partial}{\partial z_k}  \right \rangle \\ & + \sum_{j,k=1}^n \frac{2(a_j-a_k)}{a_j(a_j+a_k)}i \Big \langle \Big( (\nabla ^{X}_{\tfrac{\partial}{\partial z_j}} J) (\nabla ^{X}_{\tfrac{\partial}{\partial \bar z_j}}|\mathcal J|) + (\nabla ^{X}_{\tfrac{\partial}{\partial \bar z_j}}J)(\nabla ^{X}_{\tfrac{\partial}{\partial z_j}}|\mathcal J|) \\
& + (\nabla ^{X}_{\tfrac{\partial}{\partial \bar z_j}}|\mathcal J|) (\nabla ^{X}_{\tfrac{\partial}{\partial z_j}} J) + (\nabla ^{X}_{\tfrac{\partial}{\partial z_j}}|\mathcal J|) (\nabla ^{X}_{\tfrac{\partial}{\partial \bar z_j}}J)\Big)
\tfrac{\partial}{\partial z_k}, \tfrac{\partial}{\partial \bar z_k}  \Big \rangle  \\ 
  & -\sum_{j,k=1}^n\frac{4i}{a_j+a_k}\Big \langle  \Big((\nabla ^{X}_{\tfrac{\partial}{\partial z_j}} J)(\nabla ^{X}_{\tfrac{\partial}{\partial \bar z_k}}|\mathcal J|) + (\nabla ^{X}_{\tfrac{\partial}{\partial \bar z_k}}J) (\nabla ^{X}_{\tfrac{\partial}{\partial z_j}}|\mathcal J|) \\
& +(\nabla ^{X}_{\tfrac{\partial}{\partial \bar z_k}}|\mathcal J|) (\nabla ^{X}_{\tfrac{\partial}{\partial z_j}} J)+ (\nabla ^{X}_{\tfrac{\partial}{\partial z_j}}|\mathcal J|) (\nabla ^{X}_{\tfrac{\partial}{\partial \bar z_k}}J)\Big)
\tfrac{\partial}{\partial z_k}, \tfrac{\partial}{\partial \bar z_j} \Big \rangle\\
& + \sum_{j,k=1}^n\frac{8}{a_j+a_k}\left \langle (\nabla ^{X} \nabla ^{X}|\mathcal J|)_{(\tfrac{\partial}{\partial z_j},\tfrac{\partial}{\partial \bar z_k})}\tfrac{\partial}{\partial \bar z_j}- (\nabla ^{X} \nabla ^{X}|\mathcal J|)_{(\tfrac{\partial}{\partial z_j},\tfrac{\partial}{\partial \bar z_j})}\tfrac{\partial}{\partial \bar z_k}, \tfrac{\partial}{\partial z_k} \right \rangle \\ &  +\sum_{j,k,\ell=1}^n \frac{8}{9a_k(a_j+a_\ell)}\left|\left\langle  (\nabla ^{X}_{\tfrac{\partial}{\partial \bar z_\ell}}|\mathcal J|)\tfrac{\partial} {\partial\bar{z}_j}-(\nabla ^{X}_{\tfrac{\partial}{\partial \bar z_j}}|\mathcal J|)\tfrac{\partial} {\partial\bar{z}_\ell},\tfrac{\partial}{\partial z_k}\right\rangle\right|^2 \\ & +\sum_{j, k,\ell=1}^n \frac{2(a_\ell+a_j)^2}{a_k(a_j+a_k+a_\ell)}\left|\left\langle(\nabla ^{X}_{\tfrac{\partial}{\partial \bar z_k}}J)\tfrac{\partial}{\partial \bar z_\ell},\tfrac{\partial} {\partial\bar{z}_j}\right\rangle\right|^2.
\end{aligned}
\end{multline*}
This formula allows us to write the formula \eqref{e:rho} for the particular cases mentioned in Examples \ref{ex:1} and \ref{ex:2}.   

\section{Acknowlidgements}
This work was supported by the Laboratory of Topology and Dynamics, Novosibirsk State University (contract no. 14.Y26.31.0025 with the Ministry of Education and Science of the Russian Federation).

We are grateful to X. Ma and G. Marinescu for useful discussions.

\end{document}